\documentclass[12pt,psamsfonts,reqno]{amsart}
\usepackage{amssymb,eucal,graphics,latexsym}

\theoremstyle{plain}
\newtheorem{lemma}{Lemma}[section]
\newtheorem{theorem}[lemma]{Theorem}
\newtheorem{proposition}[lemma]{Proposition}
\newtheorem{corollary}[lemma]{Corollary}
\newtheorem{example}[lemma]{Example}
\newtheorem{claim}{Claim}

\newtheorem*{stat}{\name}
\newcommand{\name}{testing}

\theoremstyle{definition}
\newtheorem{definition}[lemma]{Definition}

\theoremstyle{remark}

\newenvironment{all}[1]{\renewcommand{\name}{#1}\begin{stat}}
                        {\end{stat}}

\newcommand{\qedc}{{\qed}~{\rm Claim~{\theclaim}.}}

\newenvironment{cproof}
{\begin{proof}[Proof of Claim.]}
{\qedc\renewcommand{\qed}{}\end{proof}}

\numberwithin{equation}{section}

\newcommand{\famm}[2]{\left({#1}\mid{#2}\right)}
\newcommand{\set}[1]{\{{#1}\}}
\newcommand{\setm}[2]{\{{#1}\mid{#2}\}}

\newcommand{\pup}[1]{\textup{(}{#1}\textup{)}}

\newcommand{\ol}[1]{\overline{{#1}}}
\newcommand{\jz}{$(\vee,0)$}

\newcommand{\jzs}{\jz-semi\-lat\-tice}

\newcommand{\mh}{meet-ho\-mo\-mor\-phism}
\newcommand{\jh}{join-ho\-mo\-mor\-phism}

\newcommand{\js}{join-sem\-i\-lat\-tice}

\newcommand{\ms}{meet-sem\-i\-lat\-tice}
\newcommand{\mss}{meet-sub\-sem\-i\-lat\-tice}
\newcommand{\res}{\mathbin{\restriction}}
\newcommand{\Pow}{\mathfrak{P}}
\DeclareMathOperator{\dom}{dom}
\DeclareMathOperator{\rng}{rng}

\DeclareMathOperator{\crit}{crit}

\newcommand{\cA}{\mathcal{A}}
\newcommand{\cB}{\mathcal{B}}
\newcommand{\cD}{\mathcal{D}}
\newcommand{\cF}{\mathcal{F}}
\newcommand{\cK}{\mathcal{K}}
\newcommand{\cM}{\mathcal{M}}
\newcommand{\cW}{\mathcal{W}}
\newcommand{\cZ}{\mathcal{Z}}

\newcommand{\ZFC}{\mathsf{ZFC}}
\newcommand{\MA}{\mathsf{MA}}
\newcommand{\MAP}{\mathsf{MA}(\aleph_1;\text{precaliber }\aleph_1)}
\newcommand{\xL}{\mathbf{L}}
\newcommand{\xV}{\mathbf{V}}
\newcommand{\id}{\mathrm{id}}
\newcommand{\es}{\varnothing}
\newcommand{\into}{\hookrightarrow}
\newcommand{\onto}{\twoheadrightarrow}
\newcommand{\les}{\leqslant}

\newcommand{\dnw}{\mathbin{\downarrow}}

\newcommand{\utr}{\trianglelefteq}
\newcommand{\tr}{\vartriangleleft}
\newcommand{\go}{\omega}
\newcommand{\GG}{\mathbb{G}}

\DeclareMathOperator{\Proj}{Proj}
\DeclareMathOperator{\Sk}{Sk}

\author[F.~Wehrung]{Friedrich Wehrung}
\address{LMNO, CNRS UMR 6139\\
D\'epartement de Math\'ematiques, BP 5186\\
Universit\'e de Caen, Campus 2\\
14032 Caen cedex\\
France}
\email{wehrung@math.unicaen.fr}
\urladdr{http://www.math.unicaen.fr/\~{}wehrung}

\date{\today}

\title{Large semilattices of breadth three}

\begin{document}
\baselineskip=17pt     %[after \begin{document}]

\begin{abstract}
A 1984 problem of S.\,Z. Ditor asks whether there exists a lattice of cardinality~$\aleph_2$, with zero, in which every principal ideal is finite and every element has at most three lower covers. We prove that the existence of such a lattice follows from either one of two axioms that are known to be independent of~$\ZFC$, namely (1)~Martin's Axiom restricted to collections of~$\aleph_1$ dense subsets in posets of precaliber~$\aleph_1$, (2)~the existence of a gap-$1$ morass. In particular, the existence of such a lattice is consistent with~$\ZFC$, while the non-existence of such a lattice implies that~$\omega_2$ is inaccessible in the constructible universe.

We also prove that for each regular uncountable cardinal~$\kappa$ and each positive integer~$n$, there exists a \jzs\ $L$ of cardinality~$\kappa^{+n}$ and breadth $n+1$ in which every principal ideal has less than~$\kappa$ elements.
\end{abstract}

\maketitle

\section{Introduction}\label{S:Intro}
Various representation theorems, stating that every object of `size'~$\aleph_1$ belongs to the range of a given functor, rely on the existence of lattices called \emph{$2$-ladders}. By definition, a $2$-ladder is a lattice with zero, in which every principal ideal is finite, and in which every element has at most two lower covers. Every $2$-ladder has cardinality at most $\aleph_1$, and the existence of $2$-ladders of cardinality exactly~$\aleph_1$ was proved in Ditor~\cite{Dito84} (cf. Proposition~\ref{P:2ladd}). These $2$-ladders have been used in various contexts such as abstract measure theory (Dobbertin~\cite{Dobb86}), lattice theory (Gr\"atzer, Lakser, and Wehrung~\cite{ConAm}), ring theory (Wehrung~\cite{Al1}), or general algebra (R\r{u}\v{z}i\v{c}ka, T\r{u}ma, and Wehrung~\cite{RTW}). A sample result, established in~\cite{RTW}, states that \emph{Every distributive algebraic lattice with at most~$\aleph_1$ compact elements is isomorphic to the lattice of all normal subgroups of some locally finite group}. (Here and in many related results, the~$\aleph_1$ bound turns out to be optimal.)

The basic idea of ``ladder proofs'' is always the same: we are given categories~$\cA$ and~$\cB$ together with a functor~$\Phi\colon\cA\to\cB$ and a `large' object~$S$ of~$\cB$ (of `size'~$\aleph_1$), that we wish to represent as~$\Phi(X)$, for some object~$X$ in the domain of~$\Phi$. We represent~$S$ as a direct limit $S=\varinjlim_{i\in I}S_i$ of (say) `finite' objects~$S_i$, where~$I$ is an upward directed poset of cardinality~$\aleph_1$ (often the lattice of all finite subsets of~$S$ in case we are dealing with a concrete category). Then, using the existence of a $2$-ladder of cardinality~$\aleph_1$, we can replace the original poset~$I$ by that $2$-ladder. Under ``amalgamation-type'' conditions, this makes it possible to represent each~$S_i$ as some~$\Phi(X_i)$, with transition morphisms between the~$X_i$s being constructed in such a way that $X=\varinjlim_{i\in I}X_i$ can be defined and $\Phi(X)\cong S$.

We define \emph{$3$-ladders} the same way as $2$-ladders, except that ``two lower covers'' is replaced by ``three lower covers''. The problem of existence of $3$-ladders of cardinality~$\aleph_2$ was posed in Ditor~\cite{Dito84}. Such $3$-ladders would presumably be used in trying to represent objects of size~$\aleph_2$. Nevertheless I must reluctantly admit that no potential use of the existence of $3$-ladders of cardinality~$\aleph_2$ had been found so far, due to the failure of a certain three-dimensional amalgamation property, of set-theoretical nature, stated in Section~10 in Wehrung~\cite{Con2}, thus making Ditor's problem quite `romantic' (and thus, somewhat paradoxically therefore arguably, attractive).

However, this situation has been evolving recently. For classes~$\cA$ and~$\cB$ of algebras, the \emph{critical point} $\crit(\cA;\cB)$ is defined in Gillibert~\cite{Gill07} as the least possible cardinality of a semilattice in the compact congruence class of~$\cA$ but not of~$\cB$ if it exists (and, say, $\infty$ otherwise). In case both~$\cA$ and~$\cB$ are finitely generated lattice varieties, it is proved in Gillibert~\cite{Gill07} that $\crit(\cA;\cB)$ is either finite, or~$\aleph_n$ for some natural number~$n$, or~$\infty$. In the second case only examples with~$n\in\set{0,1,2}$ have been found so far (Plo\v{s}\v{c}ica~\cite{Plos00,Plos03}, Gillibert~\cite{Gill07}). Investigating the possibility of~$n=3$ (i.e., $\crit(\cA;\cB)=\aleph_3$) would quite likely require $3$-ladders of cardinality~$\aleph_2$.

The present paper is intended as an encouragement in that direction. We partially solve one of Ditor's problems by giving a (rather easy) proof that for each regular uncountable cardinal~$\kappa$ and each positive integer~$n$, there exists a \jzs\ of cardinality~$\kappa^{+n}$ and breadth~$n+1$ in which every principal ideal has less than~$\kappa$ elements (cf. Theorem~\ref{T:LargeBdedBr}). Furthermore, although we are still not able to settle whether the existence of a $3$-ladder of cardinality~$\aleph_2$ is provable in the usual axiom system~$\ZFC$ of set theory with the Axiom of Choice, we prove that it follows from either one of two quite distinct, and in some sense `orthogonal', set-theoretical axioms, namely a weak form of Martin's Axiom plus $2^{\aleph_0}>\aleph_1$ denoted by~$\MAP$ (cf. Theorem~\ref{T:MAP23ladd}) and the existence of a gap-$1$ morass (cf. Theorem~\ref{T:Mor2Ladd}). In particular, the \emph{existence} of a $3$-ladder of cardinality~$\aleph_2$ is consistent with~$\ZFC$, while the \emph{non-existence} of a $3$-ladder of cardinality~$\aleph_2$ implies that~$\omega_2$ is inaccessible in the constructible universe.

Our proofs are organized in such a way that no prerequisites in lattice theory and set theory other than the basic ones are necessary to read them. Hence we hope to achieve intelligibility for both lattice-the\-o\-ret\-i\-cal and set-the\-o\-ret\-i\-cal communities.

\section{Basic concepts}\label{S:Basic}

We shall use standard set-theoretical notation and terminology. Throughout the paper, `countable' will mean `at most countable'. A cardinal is an initial ordinal, and we denote by~$\kappa^+$ the successor cardinal of a cardinal~$\kappa$. More generally, we denote by~$\kappa^{+n}$ the $n^{\mathrm{th}}$ successor cardinal of~$\kappa$, for each natural number~$n$.
We denote by~$\dom f$ (resp., $\rng f$) the domain (resp., range) of a function~$f$, and we put $f[X]:=\setm{f(x)}{x\in X}$ for each $X\subseteq\dom f$. We denote by~$\Pow(X)$ the powerset of a set~$X$, and by~$\sqcup$ the partial operation of disjoint union. A function~$f$ is \emph{finite-to-one} if the inverse image of any singleton under~$f$ is finite. For a set~$\Omega$ and a cardinal~$\lambda$, we set
 \begin{align*}
 [\Omega]^{\lambda}&=\setm{X\in\Pow(\Omega)}{|X|=\lambda}\,,\\
 [\Omega]^{\les\lambda}&=\setm{X\in\Pow(\Omega)}{|X|\leq\lambda}\,,\\
 [\Omega]^{<\lambda}&=\setm{X\in\Pow(\Omega)}{|X|<\lambda}\,.
 \end{align*}
Two elements~$x$ and~$y$ in a poset~$P$ are \emph{comparable} if either $x\leq y$ or $y\leq x$. We say that~$x$ is a \emph{lower cover} of~$y$, if $x<y$ and there is no element~$z\in P$ such that $x<z<y$; in addition, if~$x$ is the least element of~$P$ (denoted by $0_P$ if it exists), we say that~$y$ is an \emph{atom} of~$P$. We say that~$P$ is \emph{atomistic} if every element of~$P$ is a join of atoms of~$P$. For a subset~$X$ and an element~$p$ in~$P$, we set
 \[
 X\dnw p:=\setm{x\in X}{x\leq p}\,.
 \]
We say that~$X$ is a \emph{lower subset} of~$P$, if $P\dnw x\subseteq X$ for each $x\in X$. (In forcing terminology, this means that~$X$ is \emph{open}.) An \emph{ideal} of~$P$ is a nonempty, upward directed, lower subset of~$P$; it is a \emph{principal ideal} if it is equal to~$P\dnw p$ for some~$p\in P$.
A \emph{filter} of~$P$ is an ideal of the dual poset of~$P$.
We say that~$P$ is \emph{lower finite} if $P\dnw p$ is finite for each $p\in P$. 
An \emph{order-embedding} from a poset~$P$ into a poset~$Q$ is a map $f\colon P\to Q$ such that $f(x)\leq f(y)$ if{f} $x\leq y$, for all $x,y\in P$. An order-embedding~$f$ is a \emph{lower embedding} if the range of~$f$ is a lower subset of~$Q$. Observe that a lower embedding preserves all meets of nonempty subsets in~$P$, and all joins of nonempty finite subsets of~$P$ in case~$Q$ is a \js.

For an element~$x$ and a subset~$F$ in~$P$, we denote by~$x^F$ the least element of~$F$ above~$x$ if it exists.

\v{S}anin's classical \emph{$\Delta$-Lemma} (cf. Jech \cite[Theorem~9.18]{Jech03}) is the following.

\begin{all}{$\Delta$-Lemma}
Let~$\cW$ be an uncountable collection of finite sets. Then there are an uncountable subset~$\cZ$ of~$\cW$ and a finite set~$R$ \pup{the \emph{root} of~$\cZ$} such that $X\cap Y=R$ for all distinct~$X,Y\in\cZ$.
\end{all}

We recall some basic terminology in the theory of forcing, which we shall use systematically in Section~\ref{S:PPK}. A subset~$X$ in a poset~$P$ is \emph{dense} if it meets every principal ideal of~$P$. We say that~$X$ is \emph{centred} if every finite subset of~$X$ has a lower bound in~$P$. We say that~$P$ has \emph{precaliber~$\aleph_1$} if every uncountable subset of~$P$ has an uncountable centred subset; in particular, this implies the countable chain condition. For a collection~$\cD$ of subsets of~$P$, a subset~$G$ of~$P$ is \emph{$\cD$-generic} if $G\cap D\neq\es$ for each \emph{dense} $D\in\cD$. The following classical lemma (cf. Jech \cite[Lemma~14.4]{Jech03}) is, formally, a poset analogue of Baire's category Theorem.

\begin{lemma}\label{L:CtbleGen}
Let $\cD$ be a countable collection of subsets of a poset~$P$. Then each element of~$P$ is contained in a $\cD$-generic filter on~$P$.
\end{lemma}

For~$|\cD|=\aleph_1$ there may be no $\cD$-generic filters of~$P$ (cf. Jech \cite[Exercise~16.11]{Jech03}), nevertheless for restricted classes of partial orderings~$P$ we obtain set-theoretical axioms that are independent of~$\ZFC$. We shall need the following proper weakening of the~$\aleph_1$ instance of Martin's Axiom~$\MA$ usually denoted by either $\MA_{\aleph_1}$ or~$\MA(\aleph_1)$ (cf. Jech \cite[Section~16]{Jech03} and Weiss \cite[Section~3]{Weiss}); the first parameter~$\aleph_1$ refers to the cardinality of~$\cD$.

\begin{all}{$\MAP$}
For every poset~$P$ of precaliber~$\aleph_1$ and every collection~$\cD$ of subsets of~$P$, if $|\cD|\leq\aleph_1$, then there exists a $\cD$-generic filter on~$P$.
\end{all}

\section{Simplified morasses}\label{S:Morass}
For a positive integer~$n$, trying to build certain structures of size~$\aleph_{n+1}$ as direct limits of countable structures may impose very demanding constraints on the direct systems used for the construction. The pattern of the repetitions of the countable building blocks and their transition morphisms in the direct system is then coded by a complex combinatorial object called a \emph{gap-$n$ morass}. Gap-$n$ morasses were introduced by Ronald Jensen in the seventies, enabling him to solve positively the \emph{finite gap cardinal transfer conjecture} in the constructible universe~$\xL$. The existence of morasses is independent of the usual axiom system of set theory $\ZFC$. For example, there are gap-$1$ morasses in~$\xL$ (cf. Devlin \cite[Section~VIII.2]{Devl84}), and even in the universe~$\xL[A]$ of sets constructible with oracle~$A$, for any $A\subseteq\omega_1$; hence if~$\omega_2$ is not inaccessible in~$\xL$, then there is a gap-$1$ morass in the ambient set-theoretical universe~$\xV$ (cf. Devlin \cite[Exercise~VIII.6]{Devl84}). Conversely, the existence of a gap-$1$ morass implies the existence of a Kurepa tree while, in the generic extension obtained by Levy collapsing an inaccessible cardinal on~$\omega_2$ while preserving~$\omega_1$, there is no Kurepa tree (cf. Silver~\cite{Silv71}). In particular, the non-existence of a gap-$1$ morass is equiconsistent, relatively to $\ZFC$, to the existence of an inaccessible cardinal.

However, even for $n=1$ the combinatorial theorems involving morasses are hard to come by, due to the extreme complexity of the definition of gap-$n$ morasses. Fortunately, the definition of a gap-$1$ morass has been greatly simplified by Dan Velleman~\cite{Vell84a}, where it is proved that the existence of a gap-$1$ morass is equivalent to the existence of a `simplified $(\omega_1,1)$-morass'.

We denote the (noncommutative) ordinal addition by~$+$.
For ordinals $\alpha\leq\nobreak\beta$, we denote by $\beta-\alpha$ the unique ordinal~$\xi$ such that $\alpha+\xi=\beta$. Furthermore, we denote by~$\tau_{\alpha,\beta}$ the order-embedding from ordinals to ordinals defined by 
 \begin{equation}\label{Eq:Deftauab}
 \tau_{\alpha,\beta}(\xi)=\begin{cases}
 \xi\,,&\text{if }\xi<\alpha\,,\\
 \beta+(\xi-\alpha)\,,&\text{if }\xi\geq\alpha\,.
 \end{cases}
 \end{equation}
We shall use the following definition, obtained by slightly amending the one in Devlin \cite[Section~VIII.4]{Devl84} by requiring $\theta_0=2$ instead of $\theta_0=1$, but all $\delta_\alpha$s nonzero (which could not hold for $\theta_0=1$). It is easily obtained from a simplified morass as defined in Devlin \cite[Section~VIII.4]{Devl84} by adding a new zero element to each~$\theta_\alpha$ (that is, by replacing~$\theta_\alpha$ by $1+\theta_\alpha$) and replacing any $f\in\cF_{\alpha,\beta}$ by the unique zero-preserving map sending $1+\xi$ to $1+f(\xi)$, for each $\xi<\theta_\alpha$.

\begin{definition}\label{D:Gap1Mor}
Let~$\kappa$ be an infinite cardinal.
A \emph{simplified $(\kappa,1)$-morass} is a structure
 \[
 \cM=\bigl(\famm{\theta_\alpha}{\alpha\leq\kappa},
 \famm{\cF_{\alpha,\beta}}{\alpha<\beta\leq\kappa}\bigr)
 \]
satisfying the following conditions:
\begin{itemize}
\item[(P0)]\begin{itemize}
\item[(a)] $\theta_0=2$, $0<\theta_\alpha<\kappa$ for each $\alpha<\kappa$, and $\theta_\kappa=\kappa^+$.

\item[(b)] $\cF_{\alpha,\beta}$ is a set of order-embeddings from~$\theta_\alpha$ into~$\theta_\beta$, for all $\alpha<\nobreak\beta\leq\nobreak\kappa$.
\end{itemize}

\item[(P1)] $|\cF_{\alpha,\beta}|<\kappa$, for all $\alpha<\beta<\kappa$.

\item[(P2)] If $\alpha<\beta<\gamma\leq\kappa$, then $\cF_{\alpha,\gamma}=\setm{f\circ g}{f\in\cF_{\beta,\gamma}\text{ and }g\in\cF_{\alpha,\beta}}$.

\item[(P3)] For each $\alpha<\kappa$, there exists a \emph{nonzero} ordinal $\delta_\alpha<\theta_\alpha$ such that $\theta_{\alpha+1}=\theta_\alpha+(\theta_\alpha-\delta_\alpha)$ and $\cF_{\alpha,\alpha+1}=\set{\id_{\theta_\alpha},f_\alpha}$, where $f_\alpha$ denotes the restriction of $\tau_{\delta_\alpha,\theta_\alpha}$ from~$\theta_\alpha$ into~$\theta_{\alpha+1}$.

\item[(P4)] For every limit ordinal $\lambda\leq\kappa$, all $\alpha_i<\lambda$ and $f_i\in\cF_{\alpha_i,\lambda}$, for $i<2$, there exists $\alpha<\lambda$ with $\alpha_0,\alpha_1<\alpha$ together with $f'_i\in\cF_{\alpha_i,\alpha}$, for $i<2$, and $g\in\cF_{\alpha,\lambda}$ such that $f_i=g\circ f'_i$ for each $i<2$.

\item[(P5)] The equality
$\theta_\alpha=\bigcup\famm{f[\theta_\xi]}
{\xi<\alpha\text{ and }f\in\cF_{\xi,\alpha}}$ holds for each $\alpha>\nobreak0$.
\end{itemize}
\end{definition}

It is proved in Velleman~\cite{Vell84a} that for~$\kappa$ regular uncountable, there exists a $(\kappa,1)$-morass if{f} there exists a simplified $(\kappa,1)$-morass. For the countable case, the existence of a $(\omega,1)$-morass is provable in~$\ZFC$, see Velleman~\cite{Vell84b}.

Simplified morasses as above satisfy the following simple but very useful lemma, which is the basis of the construction of the Kurepa tree obtained from a $(\kappa,1)$-morass (cf. Velleman \cite[Lemma~3.2]{Vell84a}).

\begin{lemma}\label{L:Mor2Tree}
Let~$\alpha<\beta\leq\kappa$, let $\xi_0,\xi_1<\theta_\alpha$, and let $f_0,f_1\in\cF_{\alpha,\beta}$. If $f_0(\xi_0)=f_1(\xi_1)$, then $\xi_0=\xi_1$ and $f_0\res_{\xi_0}=f_1\res_{\xi_1}$.
\end{lemma}

\section{Ladders and breadth}\label{S:Breadth}
The classical definition of breadth, see Ditor \cite[Section~4]{Dito84}, runs as follows. Let~$n$ be a positive integer. A \js\ $S$ has \emph{breadth at most~$n$} if for every nonempty finite subset~$X$ of~$L$, there exists a nonempty $Y\subseteq X$ with at most~$n$ elements such that $\bigvee X=\bigvee Y$. This is a particular case of the following definition of breadth, valid for every poset, and, in addition, self-dual: we say that a poset~$P$ has breadth at most~$n$, if for all $x_i$, $y_i$ \pup{$0\leq i\leq n$} in $L$, if $x_i\leq y_j$ for all $i\neq j$ in $\set{0,1,\dots,n}$, then there exists $i\in\set{0,1,\dots,n}$ such that $x_i\leq y_i$.

\begin{definition}\label{S:ladder}
Let $k$ be a positive integer. A \emph{$k$-ladder} is a lower finite lattice in which every element has at most~$k$ lower covers.
\end{definition}

Every $k$-ladder has breadth at most~$k$. The diamond~$M_3$ has breadth~$2$ but it is not a $2$-ladder. Every finite chain is a $1$-ladder. The chain $\go$ of all non-negative integers is also a $1$-ladder. Note that $k$-ladders are called \emph{$k$-frames} in Dobbertin~\cite{Dobb86}; the latter terminology being already used for a completely different lattice-theoretical concept (von~Neumann frames), we will not use it. The following is proved in Ditor~\cite{Dito84}.

\begin{proposition}\label{P:UppBdBr}
Let~$k$ be a positive integer. Then every lower finite lattice of breadth at most~$k$ \pup{thus, in particular, every $k$-ladder} has at most~$\aleph_{k-1}$ elements.
\end{proposition}

Proposition~\ref{P:UppBdBr} is especially easy to prove by using Kuratowski's Free Set Theorem, see Kuratowski~\cite{Kura51}.
The converse is obviously true for $k=1$---that is, there exists a $1$-ladder of cardinality~$\aleph_0$ (namely, the chain~$\omega$ of all natural numbers); also for $k=2$, by the following result of Ditor~\cite{Dito84}, also proved by Dobbertin~\cite{Dobb86}. We include a proof for convenience.

\begin{proposition}\label{P:2ladd}
There exists an atomistic $2$-ladder of cardinality $\aleph_1$.
\end{proposition}

\begin{proof}
We construct inductively a $\omega_1$-sequence $\cF=\famm{F_\alpha}{\alpha<\omega_1}$ of countable atomistic $2$-ladders, such that $\alpha<\beta$ implies that~$F_\alpha$ is a proper ideal of~$F_\beta$. Once this is done, the $2$-ladder $F_{\omega_1}:=\bigcup\famm{F_\alpha}{\alpha<\omega_1}$ will clearly solve our problem.

We take $F_0:=\set{0}$. If~$\lambda<\omega_1$ is a limit ordinal and all $F_\alpha$, for $\alpha<\lambda$, are constructed, take $F_\lambda:=\bigcup\famm{F_\alpha}{\alpha<\lambda}$. Suppose that~$F_\alpha$ is constructed. If~$F_\alpha$ is finite, pick outside objects~$p$, $1$ and put $F_{\alpha+1}:=F_\alpha\cup\set{p,1}$, with the additional relations $p<1$ and $x<1$ for each $x\in F_\alpha$. If~$F_\alpha$ is infinite, then, as it is a countable lattice, it has a strictly increasing cofinal sequence~$\famm{c_n}{1\leq n<\omega}$. Consider a one-to-one sequence $\famm{d_n}{n<\omega}$ of objects outside~$F_\alpha$ and put
 \[
 F_{\alpha+1}:=F_\alpha\cup\setm{d_n}{n<\omega}\,,
 \]
with the additional relations $c_n<d_n$ for $1\leq n<\omega$ and $d_n<d_{n+1}$ for $0\leq n<\omega$. The $\omega_1$-sequence~$\cF$ thus constructed is as required.
\end{proof}

The most natural attempt at proving the existence of a $3$-ladder of cardinality~$\aleph_2$, by imitating the proof of Proposition~\ref{P:2ladd}, would require that every $3$-ladder of cardinality~$\aleph_1$ has a cofinal \mss\ which is also a $2$-ladder (cf. the proof of Theorem~\ref{T:MAP23ladd}). We do not know whether this statement is a theorem of~$\ZFC$, although, by Theorem~\ref{T:ExGenaleph1}, it is consistent with~$\ZFC$. The following example shows that the most straightforward attempt at proving that statement, by expressing structures of cardinality~$\aleph_1$ as directed unions of countable structures, fails.

\begin{example}\label{Ex:LL'FnoF'}
There exists a countable $3$-ladder~$K'$ with an ideal~$K$ and a cofinal \mss~$F$ of~$K$ which is also a $2$-ladder, although there is no cofinal \mss\ of breadth at most two of~$K'$ containing~$F$.
\end{example}

\begin{proof}
We denote by~$K'$ the lattice represented in Figure~\ref{Fig:no2ladd}, and we put
 \begin{align*}
 K&:=K'\setminus\setm{t_n}{n<\omega}\,,\\
 F&:=\setm{x_n}{n<\omega}\cup\setm{y_n}{n<\omega}
 \cup\setm{x_n\wedge y_n}{n<\omega}\,. 
 \end{align*}
 
\begin{figure}[htb]
\includegraphics{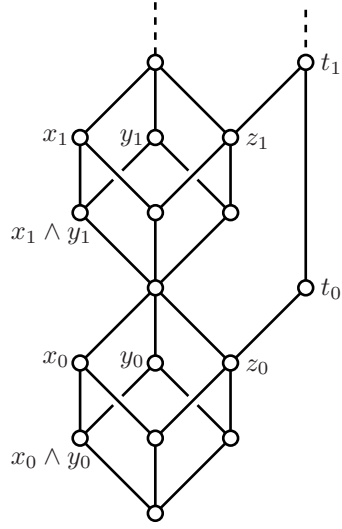}
\caption{The lattice $K'$}\label{Fig:no2ladd}
\end{figure}

Suppose that there is a cofinal \mss~$F'$ of~$K'$, with breadth at most two, containing~$F$. As~$F'$ is cofinal in~$K'$, there exists $n<\omega$ such that $t_n\in F'$. As~$x_{n+1}$ belongs to~$F$, it also belongs to~$F'$, thus ~$z_n=t_n\wedge x_{n+1}$ also belongs to~$F'$, so $\set{x_n,y_n,z_n}$ is contained in~$F'$, and so, as~$F'$ has breadth at most two, $x_n\wedge y_n\wedge z_n$ belongs to $\set{x_n\wedge y_n,x_n\wedge z_n,y_n\wedge z_n}$, a contradiction.
\end{proof}

A particular case of the problem above is stated on top of Page~58 in Ditor~\cite{Dito84}: let~$F$ be a $2$-ladder of cardinality~$\aleph_1$. Does the $3$-ladder $F\times\omega$, endowed with the product order, have a cofinal \mss\ which is also a $2$-ladder? The answer to that question is affirmative, due to the following easy result.

\begin{proposition}\label{P:Ktimesom}
Let $K$ be an infinite, lower finite lattice. Then $K\times\omega$ has a cofinal \mss\ isomorphic to~$K$.
\end{proposition}

\begin{proof}
As~$K$ is an infinite, lower finite lattice, it has a strictly increasing sequence $\famm{a_n}{n<\omega}$. As~$K$ is lower finite, for each $x\in K$, there exists a largest natural number~$n$ such that $a_n\leq x$; denote this integer by~$f(x)$. We define~$K'$ as the graph of~$f$, that is, $K':=\setm{(x,f(x))}{x\in K}$. As~$f$ is a \mh, $K'$ is a \mss\ of $K\times\omega$. Furthermore, $x\mapsto(x,f(x))$ defines a lattice isomorphism of~$K$ onto~$K'$, and~$K'$ is obviously cofinal in~$K\times\omega$.
\end{proof}

\section{Large semilattices with bounded breadth}\label{S:LargeBdedBr}

In this section we shall see that weak analogues of $k$-ladders, obtained by replacing~$\aleph_0$ by a regular uncountable cardinal and the condition that every element has at most~$k$ lower covers by the breadth being below~$k$, can be easily constructed (cf. Theorem~\ref{T:LargeBdedBr}).

The following Lemma~\ref{L:BreadthGen} is a slightly more general version of \cite[Proposition~4.3]{Dito84}. It says that breadth can be verified on the generators of a \js, the difference with~\cite[Proposition~4.3]{Dito84} lying in the definition of a generating subset.

A subset~$G$ in a join-semilattice~$L$ \emph{generates}~$L$, if every element of~$L$ is a (not necessarily finite) join of elements of~$G$. Equivalently, for all $a,b\in L$ such that $a\nleq b$, there exists $g\in G$ such that $g\leq a$ and $g\nleq b$.

\begin{lemma}\label{L:BreadthGen}
Let $L$ be a join-semilattice, let $G$ be a generating subset of $L$, and let $n$ be a positive integer. Then $L$ has breadth at most~$n$ if{f} for every subset $U\in[G]^{n+1}$, there exists $u\in U$ such that $u\leq\bigvee(U\setminus\set{u})$.
\end{lemma}

\begin{proof}
We prove the nontrivial direction. Let $A\in[L]^{n+1}$, and suppose that
$a\nleq\bigvee(A\setminus\set{a})$ for each $a\in A$. As~$G$ generates~$L$, there exists $u_a\in G$ such that $u_a\leq a$ and $u_a\nleq\bigvee(A\setminus\set{a})$, for all $a\in A$. Hence $u_a\nleq\bigvee\famm{u_x}{x\in A\setminus\set{a}}$, for each $a\in A$, a contradiction as all elements~$u_x$ belong to~$G$ and by assumption.
\end{proof}

An \emph{algebraic closure operator} on a poset~$P$ is a map $f\colon P\to P$ such that $f\circ f=f$ (we say that~$f$ is \emph{idempotent}), $x\leq f(x)$ for each $x\in P$, $x\leq y$ implies that $f(x)\leq f(y)$ for all $x,y\in P$, and for each nonempty upward directed subset~$X$ of~$P$ admitting a join, the join of~$f[X]$ exists and is equal to $f\bigl(\bigvee X\bigr)$. (We shall refer to the latter property as the \emph{join-continuity} of~$f$.)

\begin{lemma}\label{L:LargeAlCl}
For each regular uncountable cardinal~$\kappa$ and each positive integer~$n$, there exists an algebraic closure operator~$f$ on~$[\kappa^{+n}]^{<\kappa}$ such that
 \begin{gather}
 (\forall U\in[\kappa^{+n}]^{n+2})(\exists\xi\in U)\bigl(\xi\in f(U\setminus\set{\xi})\bigr)\,,
 \label{Eq:NoFree}\\
 f(X)=X\quad\text{for each }X\in[\kappa^{+n}]^{\les n}\,.\label{Eq:f(X)=X}
 \end{gather}
\end{lemma}

\begin{proof}
It follows from Kuratowski~\cite{Kura51} (see also Theorem~45.7 in Erd\H os \emph{et al.}~\cite{EHMR}) that there exists a map $f_0\colon[\kappa^{+n}]^{n+1}\to[\kappa^{+n}]^{<\kappa}$ such that
 \begin{equation}\label{Eq:f0NoFree}
 (\forall U\in[\kappa^{+n}]^{n+2})(\exists\xi\in U)
 \bigl(\xi\in f_0(U\setminus\set{\xi})\bigr)\,.
 \end{equation}
We set
 \[
 g(X)=X\cup\bigcup\famm{f_0(Y)}{Y\in[X]^{n+1}}\,,\quad
 \text{for each }X\in[\kappa^{+n}]^{<\kappa}\,.
 \]
As~$\kappa$ is regular, $g$ is a self-map of~$[\kappa^{+n}]^{<\kappa}$. It obviously satisfies both~\eqref{Eq:NoFree} and~\eqref{Eq:f(X)=X}, together with all properties defining an algebraic closure operator except idempotence. Now we set
 \[
 f(X)=\bigcup\famm{g^k(X)}{k<\omega}\,,\quad\text{for each }
 X\in[\kappa^{+n}]^{<\kappa}\,.
 \]
As~$\kappa$ is regular uncountable, $f$ is a self-map of~$[\kappa^{+n}]^{<\kappa}$. It obviously satisfies both~\eqref{Eq:NoFree} and~\eqref{Eq:f(X)=X}, together with all properties defining an algebraic closure operator except idempotence. Furthermore, for every $X\in[\kappa^{+n}]^{<\kappa}$,
 \begin{align*}
 f\circ f(X)&=\bigcup\famm{f\circ g^l(X)}{l<\omega}&&
 (\text{by the join-continuity of }f)\\
 &=\bigcup\famm{g^k\circ g^l(X)}{k,l<\omega}
 &&(\text{by the definition of }f)\\
 &=f(X)&&(\text{because }g^k\circ g^l=g^{k+l})\,, 
 \end{align*}
so $f\circ f=f$.
\end{proof}

Therefore, we obtain a positive answer to the specialization of Question~A, Page~57 in Ditor~\cite{Dito84} to regular uncountable cardinals.

\begin{theorem}\label{T:LargeBdedBr}
For each regular uncountable cardinal~$\kappa$ and each positive integer~$n$, there exists an atomistic \jzs\ $L$ of breadth~$n+1$ and cardinality~$\kappa^{+n}$ such that $|L\dnw x|<\kappa$ for each~$x\in L$.
\end{theorem}

\begin{proof}
Let~$f$ be an algebraic closure operator as in Lemma~\ref{L:LargeAlCl}. We endow $L:=\setm{f(X)}{X\in[\kappa^{+n}]^{<\omega}}$ with containment. Obviously, $L$ is a \jzs\ and $|L\dnw x|<\kappa$ for each~$x\in L$. As $|L|=\kappa^{+n}$, it follows from \cite[Theorem~5.2]{Dito84} that~$L$ has breadth at least~$n+1$.  As~$n$ is nonzero and by~\eqref{Eq:f(X)=X}, every singleton~$\set{\xi}$, for $\xi<\kappa^{+n}$, belongs to~$L$, so~$L$ is atomistic, and so $G:=\setm{\set{\xi}}{\xi<\kappa^{+n}}$ generates~$L$ in the sense required by the statement of Lemma~\ref{L:BreadthGen}. Therefore, it follows from~\eqref{Eq:NoFree} together with Lemma~\ref{L:BreadthGen} that~$L$ has breadth at most~$n+1$.
\end{proof}

For $n:=1$, it it not known whether Theorem~\ref{T:LargeBdedBr} extends to singular cardinals, for instance $\kappa:=\aleph_\omega$ (cf. Problem~2 in Ditor~\cite{Dito84}). For $n:=2$, it is not known whether Theorem~\ref{T:LargeBdedBr} extends to $\kappa:=\aleph_0$ (cf. Problem~1 in Ditor~\cite{Dito84}), although we shall prove in two different ways that a positive answer is consistent with~$\ZFC$ (cf. Theorems~\ref{T:MAP23ladd} and~\ref{T:Mor2Ladd}).

\section{Preskeletons in normed lattices}\label{L:NmdLatt}

\begin{definition}\label{D:NmdLatt}
A \emph{normed lattice} is a pair $(K,\partial)$, where $K$ is a lattice and~$\partial$ is a \jh\ from~$K$ to the ordinals (the `norm'). In such a case we put
 \begin{align*}
 K_{\xi}&:=\setm{x\in K}{\partial x=\xi}\,,\\
 K_{\les\xi}&:=\setm{x\in K}{\partial x\leq\xi}\,.
 \end{align*}
Observe that each $K_{\les\xi}$ is either empty or an ideal of~$K$ (we will say that it is an \emph{extended ideal} of~$K$). We call the subsets~$K_\xi$ the \emph{levels} of~$K$. We say that $(K,\partial)$ is \emph{transitive} if its \emph{range}, that is, the range of~$\partial$, is an ordinal. Of course, in such a case, $\partial 0=0$.
\end{definition}

Observe that conversely, every increasing well-ordered sequence $\famm{K_{\les\xi}}{\xi<\theta}$ of extended ideals of~$K$ with union~$K$ defines a norm~$\partial$, \emph{via} the rule
 \[
 \partial x:=\text{least }\xi<\theta\text{ such that }x\in K_{\les\xi}\,,
 \quad\text{for each }x\in K\,.
 \]
In case $K$ is \emph{lower finite}, the set $\setm{y\in K\dnw x}{\partial y\leq\xi}$ is a finite ideal of~$K$, for every $x\in K$ and every ordinal~$\xi$, hence it has a largest element, that we shall denote by $x_{(\xi)}$. The assignment $(x,\xi)\mapsto x_{(\xi)}$ is \emph{isotone}, that is, $x\leq y$ and $\xi\leq\eta$ implies that $x_{(\xi)}\leq y_{(\eta)}$.
We shall put
 \[
 \Proj(x)=\setm{x_{(\xi)}}{\xi\in\rng\partial}\,,\quad
 \text{for each }x\in K\,.
 \]
Observe that~$\Proj(x)$ is a \emph{chain}. Furthermore, as~$\Proj(x)$ is a subset of~$K\dnw x$, it is \emph{finite}. The binary relation~$\utr$ defined by the rule
 \begin{equation}\label{Eq:defnutr}
 x\utr y\ \Leftrightarrow\ x\in\Proj(y)\,,\quad\text{for all }x,y\in K\,,
 \end{equation}
is a partial ordering on~$K$ in which all principal ideals are finite chains.
We shall always denote by~$\partial$ the norm function on a normed lattice.

We shall repeatedly use the following easy observation.

\begin{lemma}\label{L:proj2meet}
The following implications hold, for any elements~$x$ and~$y$ in a lower finite normed lattice~$K$:
\begin{enumerate}
\item $x\leq y$ implies that $\partial x=\partial(y_{(\partial x)})$.

\item $x_{(\partial y)}\leq y$ if{f} $x\wedge y=x_{(\partial y)}$.

\end{enumerate}
\end{lemma}

\begin{proof}
(i). {}From $x\leq y$ it follows that $x\leq y_{(\partial x)}$, hence $\partial x\leq\partial(y_{(\partial x)})\leq\partial x$.

(ii). We need to prove only the direct implication. It is trivial that $x_{(\partial y)}\leq x\wedge y$. For the converse, observe that $\partial(x\wedge y)\leq\partial y$, thus, as $x\wedge y\leq x$, we obtain that $x\wedge y\leq x_{(\partial y)}$, and so the equality holds.
\end{proof}

\begin{definition}\label{D:Norm2Ladd}
A \emph{preskeleton} of a lower finite normed lattice~$K$ is a subset~$F$ of~$K$ satisfying the following conditions:
\begin{enumerate}
\item $F\cap K_\xi$ is a (possibly empty) chain, for every ordinal~$\xi$;

\item $F$ is \emph{projectable}, that is, $x_{(\xi)}$ belongs to~$F$, for any $x\in F$ and any ordinal~$\xi$.
\end{enumerate}

If, in addition, $F\cap K_\xi$ is cofinal in~$K_\xi$ for each~$\xi$, we say that~$F$ is a \emph{skeleton of~$K$}.
\end{definition}

\begin{lemma}\label{L:Norm2Ladd}
Every preskeleton~$F$ of a lower finite normed lattice~$K$ satisfies the following properties:
\begin{enumerate}
\item $F$ is a \mss\ of~$K$;

\item every element of $F$ has at most two lower covers in~$F$.
\end{enumerate}
In particular, if~$F$ is upward directed, then it is a $2$-ladder.
\end{lemma}

\begin{proof}
(i). Let $x,y\in F$ and put $\alpha:=\partial(x\wedge y)$. Hence $x\wedge y=x_{(\alpha)}\wedge y_{(\alpha)}$. {}From $x\wedge y\leq x_{(\alpha)}$ it follows that $\alpha=\partial(x\wedge y)\leq\partial(x_{(\alpha)})$, thus $\partial(x_{(\alpha)})=\alpha$. Similarly, $\partial(y_{(\alpha)})=\alpha$. As both $x_{(\alpha)}$ and $y_{(\alpha)}$ belong to~$F$, it follows that they are comparable, and therefore $x\wedge y\in\set{x_{(\alpha)},y_{(\alpha)}}\subseteq F$.

(ii). Let $x\in F$. We denote by $x_*$ the largest element of~$F\cap K_{\partial x}$ smaller than~$x$ if it exists (i.e., as $F\cap K_{\partial x}$ is a chain, if $F\cap K_{\partial x}$ has an element smaller than~$x$), and by $x^-$ the largest element of~$\Proj(x)\setminus\set{x}$ if it exists (i.e., if there exists an element in~$\Proj(x)$ smaller than~$x$). It suffices to prove that every element $y\in F$ smaller than~$x$ lies either below~$x_*$ or below~$x^-$. {}From $y\leq x$ it follows that $\partial y\leq\partial x$ and $y\leq x_{(\partial y)}$. If $\partial y=\partial x$ then, as $y<x$, $x_*$ exists and $y\leq x_*$. If $\partial y<\partial x$, then, as $y\leq x$ and by Lemma~\ref{L:proj2meet}(i), $\partial(x_{(\partial y)})=\partial y<\partial x$, thus~$x^-$ exists and~$y\leq x_{(\partial y)}\leq x^-$.

Now let~$F$ be upward directed. Together with~(i) and the lower finiteness of~$F$, this implies that~$F$ is a lattice. Therefore, by~(ii), $F$ is a $2$-ladder.
\end{proof}

In particular, every skeleton of~$K$ is both a $2$-ladder and a meet-sub\-semi\-lat\-tice of~$K$. Easy examples show that a skeleton of~$K$ may not be a \emph{join}-sub\-semi\-lat\-tice of~$K$.

\section{The poset~$\Sk(K)$ of all finite preskeletons of~$K$}\label{S:PPK}

As the present section involves Martin's Axiom in an essential way, we shall use forcing terminology (\emph{open}, \emph{dense}) rather than poset terminology (\emph{lower subset}, \emph{coinitial}) throughout (cf. Section~\ref{S:Basic}).

\begin{definition}\label{D:PPK}
Let $K$ be a lower finite normed lattice (cf. Definition~\ref{D:NmdLatt}). We denote by~$\Sk(K)$ the set of all \emph{finite} preskeletons of~$K$ (cf. Definition~\ref{D:Norm2Ladd}), ordered under \emph{reverse} containment. Furthermore, we put
 \[
 \Sk_a(K):=\setm{F\in\Sk(K)}{(\exists x\in F)
 (a\leq x\text{ and }\partial a=\partial x)},\quad\text{for each }a\in K\,,
 \]
and $\cD_K=\setm{\Sk_a(K)}{a\in K}$.
\end{definition}

It is clear that~$\Sk_a(K)$ is an open subset of~$\Sk(K)$.

\begin{lemma}\label{L:PPKadense}
The subset~$\Sk_a(K)$ is dense in~$\Sk(K)$, for each~$a\in K$.
\end{lemma}

\begin{proof}
Let $E\in\Sk(K)$, we must find~$F\in\Sk_a(K)$ containing~$E$. Put\linebreak $b:=a\vee\bigvee E$. As the finite chain~$\Proj(b)$ is projectable, the subset $F:=E\cup\Proj(b)$ is projectable. Furthermore, for each $x\in E$, from $x\leq b$ it follows that $x\leq b_{(\partial x)}$, and so~$F$ is a preskeleton of~$F$.

Put $\alpha:=\partial a$. {}From $a\leq b$ it follows that $a\leq b_{(\alpha)}$ and $\partial b_{(\alpha)}=\alpha$ (cf. Lemma~\ref{L:proj2meet}(i)). Therefore, the element~$b_{(\alpha)}$ witnesses that~$F$ belongs to~$\Sk_a(K)$.
\end{proof}

\begin{definition}\label{D:LocCtble}
A normed lower finite lattice is \emph{locally countable} if all its levels are countable.
\end{definition}

\begin{proposition}\label{L:LocCtble}
Every locally countable normed lower finite lattice $K$ has cardinality at most~$\aleph_1$ and range of order-type at most~$\omega_1$.
\end{proposition}

\begin{proof}
Without loss of generality, $\theta:=\rng\partial$ is an ordinal. For each~$\alpha<\theta$ and each $a\in K_\alpha$, the map ($K_{\les\alpha}\to K_\alpha$, $x\mapsto x\vee a$) is finite-to-one, thus, as~$K_\alpha$ is countable, $K_{\les\alpha}$ is countable as well. As $|K_{\les\alpha}|\geq|\alpha|$, it follows that $\theta\leq\omega_1$. Therefore, $K=\bigcup\famm{K_\alpha}{\alpha<\theta}$ has cardinality at most~$\aleph_1$.
\end{proof}

In the following lemma, we shall use the $\Delta$-Lemma together with the notions ``centred'' and ``precaliber'' (cf. Section~\ref{S:Basic}).

\begin{lemma}\label{L:PrecalPPK}
Suppose that~$K$ is locally countable. Then the poset~$\Sk(K)$ has precaliber~$\aleph_1$.
\end{lemma}

\begin{proof}
Let $\famm{F_\alpha}{\alpha<\omega_1}$ be an $\omega_1$-sequence of elements of~$\Sk(K)$, we must find an uncountable $U\subseteq\omega_1$ such that $\setm{F_\alpha}{\alpha\in U}$ is centred. Put $e_\alpha:=\bigvee F_\alpha$ (an element of~$K$), $X_\alpha:=K\dnw e_\alpha$ (a finite subset of~$K$), and $S_\alpha=\partial[F_\alpha]$ (a finite subset of~$\rng\partial$), for each~$\alpha<\omega_1$. Two successive applications of the $\Delta$-Lemma yield an uncountable subset~$U_1$ of~$\omega_1$ and finite sets~$X$, $S$ such that
 \begin{equation}\label{Eq:FirstDelta}
 (\forall\alpha\neq\beta\text{ in }U_1)(X_\alpha\cap X_\beta=X\text{ and }
 S_\alpha\cap S_\beta=S)\,.
 \end{equation}
As~$S$ is a finite subset of~$\rng\partial$, it is contained in~$\delta+1$ for some $\delta\in\rng\partial$. For each $\alpha\in U_1$, $F_\alpha\cap K_{\les\delta}$ is a finite subset of~$K_{\les\delta}$, hence, as $K_{\les\delta}$ is countable, there are a finite subset~$F$ of~$K_{\les\delta}$ and an uncountable subset~$U_2$ of~$U_1$ such that
 \begin{equation}\label{Eq:DefnF}
 (\forall\alpha\in U_2)(F_\alpha\cap K_{\les\delta}=F)\,.
 \end{equation}
We claim that $\setm{F_\alpha}{\alpha\in U_2}$ is centred. It is sufficient to prove that for each nonempty finite subset~$A$ of~$U_2$, the union $\ol{F}:=\bigcup\famm{F_\alpha}{\alpha\in A}$ belongs to~$\Sk(K)$. As~$\ol{F}$ is obviously projectable, it suffices to prove that any elements $x,y\in\ol{F}\cap K_\xi$, for $\xi\in\rng\partial$, are comparable. Let $\alpha,\beta\in A$ such that $(x,y)\in F_\alpha\times F_\beta$. If $\alpha=\beta$ then, as~$F_\alpha$ is a preskeleton of~$K$, we are done. Suppose that $\alpha\neq\beta$. As $\xi=\partial x=\partial y$ belongs to $S_\alpha\cap S_\beta=S$ (cf. \eqref{Eq:FirstDelta}), we get $\xi\leq\delta$, thus, by~\eqref{Eq:DefnF}, $x\in F_\alpha\cap K_{\les\delta}=F_\beta\cap K_{\les\delta}\subseteq F_\beta$, so $\set{x,y}\subseteq F_\beta\cap K_\xi$, and so, as $F_\beta\cap K_\xi$ is a chain, $x$ and~$y$ are comparable.
\end{proof}

The following lemma describes $\cD_K$-generic filters over~$\Sk(K)$.

\begin{lemma}\label{L:GenonPPK}
Let $K$ be a lower finite normed lattice and let $\GG$ be a~$\cD_K$-generic filter of~$\Sk(K)$. Then $G:=\bigcup\GG$ is a skeleton of~$K$. In particular, $G$ is a cofinal \ms\ in~$K$ and it is a $2$-ladder.
\end{lemma}

\begin{proof}
As~$\GG$ is an upward directed (for containment) set of preskeletons of~$K$, $G$ is also a preskeleton of~$K$. As~$\GG$ meets~$\Sk_a(K)$ for each $a\in K$, $G$ is cofinal on each level of~$K$. In particular, $G$ is cofinal in~$K$, thus it is upward directed, and thus (cf. Lemma~\ref{L:Norm2Ladd}) it is a $2$-ladder.
\end{proof}

An immediate application of Lemmas~\ref{L:CtbleGen}, \ref{L:PPKadense}, and \ref{L:PrecalPPK} yields the following theorem. However, as there is an easy direct proof, we provide it as well.

\begin{theorem}\label{T:large2ladd}
Every countable lower finite normed lattice~$K$ has a skeleton.
\end{theorem}

\begin{proof}
As~$K$ is countable, it has a cofinal chain~$C$. It is obvious that the subset $F:=\bigcup\famm{\Proj(c)}{c\in C}$ is projectable. For any~$\xi\in\rng\partial$, any two elements of~$F\cap K_\xi$ have the form~$c_{(\xi)}$ and~$d_{(\xi)}$, for some $c,d\in C$. As~$c$ and~$d$ are comparable, so are~$c_{(\xi)}$ and~$d_{(\xi)}$. This proves that~$F$ is a preskeleton of~$K$. As~$C$ is cofinal in~$K$, every~$x\in K$ lies below some~$c\in C$. Hence $x\leq c_{(\partial x)}$ and, by Lemma~\ref{L:proj2meet}(i), $\partial x=\partial(c_{(\partial x)})$. Therefore, $F\cap K_{\partial x}$ is cofinal in~$K_{\partial x}$.
\end{proof}

For lattices of cardinality~$\aleph_1$ we get the following result:

\begin{theorem}\label{T:ExGenaleph1}
Assume that the axiom $\MAP$ holds.
Let $K$ be a lower finite lattice of cardinality at most~$\aleph_1$. Then~$K$ has a skeleton.
\end{theorem}

\begin{proof}
Write $K=\setm{e_\xi}{\xi<\omega_1}$, and denote by $I_\alpha$ the ideal of~$K$ generated by $\setm{e_\xi}{\xi<\alpha}$, for each $\alpha<\omega_1$. Finally, for each $x\in K$, denote by $\partial x$ the least ordinal~$\alpha$ such that $x\in I_\alpha$. Then $(K,\partial)$ is a locally countable, lower finite normed lattice. Now apply Lemmas~\ref{L:PPKadense} and \ref{L:PrecalPPK}.
\end{proof}

We obtain the following result.

\begin{theorem}\label{T:MAP23ladd}
Suppose that $\MAP$ holds. Then there exists an atomistic $3$-ladder of cardinality~$\aleph_2$.
\end{theorem}

\begin{proof}
We argue as in the proof of Proposition~\ref{P:2ladd}.
We construct inductively a $\omega_2$-sequence $\cK=\famm{K_\alpha}{\alpha<\omega_2}$ of atomistic $3$-ladders of cardinality at most~$\aleph_1$, such that $\alpha<\beta$ implies that~$K_\alpha$ is a proper ideal of~$K_\beta$. Once this is done, the $3$-ladder $K_{\omega_2}:=\bigcup\famm{K_\alpha}{\alpha<\omega_2}$ will clearly solve our problem.

We take $K_0:=\set{0}$. If~$\lambda<\omega_2$ is a limit ordinal and all $K_\alpha$, for $\alpha<\lambda$, are constructed, take $K_\lambda:=\bigcup\famm{K_\alpha}{\alpha<\lambda}$. Let $\alpha<\omega_2$ and suppose that~$K_\alpha$ is constructed, of cardinality at most~$\aleph_1$. By Theorem~\ref{T:ExGenaleph1} and Lemma~\ref{L:Norm2Ladd}, $K_\alpha$ has a cofinal \mss~$F$ which is also a $2$-ladder. We consider an isomorphic copy~$F^*$ of~$F$ disjoint (set-theoretically) from~$K_\alpha$, \emph{via} an isomorphism $f\colon F\onto F^*$, and we consider the partial ordering on~$K_{\alpha+1}:=K_\alpha\sqcup F^*$ obtained as the union of the respective partial orderings of~$K_\alpha$ and $F^*$ together with the additional pairs
 \[
 x<f(y)\text{ just in case }x\leq y\,,\quad
 \text{for all }(x,y)\in K_\alpha\times F\,.
 \]
As~$F$ is a \mss\ of~$K_\alpha$, $K_{\alpha+1}$ is a lattice.
It is easily seen to be an atomistic $3$-ladder (the only atom in $K_{\alpha+1}$ not in~$K_\alpha$ is $f(0_F)$). In particular, the lower covers of $f(x)$, for $x\in F$, are $x$ and (in case $x>0_F$) $f(x_0)$, $f(x_1)$, where ~$x_0$ and~$x_1$ are the lower covers of~$x$ in~$F$.

The $\omega_2$-sequence~$\cK$ thus constructed is as required.
\end{proof}

\section{Amalgamating countable normed $3$-ladders}\label{S:CtbleNdLatt}

\begin{definition}\label{D:StrAmalg}
A poset $P$ is the \emph{strong amalgam} of two subsets~$A$ and~$B$ over a subset~$I$, if $P=A\cup B$, $I=A\cap B$, and for all $(a,b)\in A\times B$, $a\leq b$ (resp., $a\geq b$) if{f} there exists $x\in I$ such that $a\leq x\leq b$ (resp., $a\geq x\geq b$).
\end{definition}

The proof of the existence of a $3$-ladder of cardinality~$\aleph_2$ from a $(\omega_1,1)$-morass is based on the following lemma (the maps~$\tau_{\delta,\theta}$ are defined in~\eqref{Eq:Deftauab}).

\begin{lemma}\label{L:amalg}
Let $K$ be a countable, transitive, normed $3$-ladder with range~$\theta$, let $\delta$ be an ordinal with $0<\delta<\theta$, and put $I:=\setm{x\in K}{\partial x<\delta}$. Then there are a countable, transitive, normed $3$-ladder~$\ol{K}$ with range $\theta':=\theta+(\theta-\delta)$, containing~$K$ as an ideal, and a lower embedding~$f\colon K\into\ol{K}$ such that the following conditions hold:
\begin{enumerate}
\item $K\cup f[K]$ is the strong amalgam of~$K$ and $f[K]$ over~$I$;

\item every element of $\ol{K}$ is the join of two elements of $K\cup f[K]$;

\item $f\res_I$ is the identity map on~$I$;

\item The equality $\partial f(x)=\tau_{\delta,\theta}(\partial x)$ holds for each $x\in K$.
\end{enumerate}
\end{lemma}

\begin{proof}
It is straightforward to construct a (set-theoretical) copy~$K^*$ of~$K$ such that $K\cap K^*=I$, isomorphic to~$K$ \emph{via} a bijection $f\colon K\onto K^*$ such that $f\res_I=\id_I$. Defining the ordering on $K\cup K^*$ according to Definition~\ref{D:StrAmalg}, we obtain that $K\cup K^*$ is a strong amalgam of~$K$ and $K^*=f[K]$ over~$I$. (As~$\delta>0$, $0_K=f(0_K)$ remains the least element of~$K\cup K^*$.) We extend the norm~$\partial$ to~$K\cup K^*$ by setting
 \begin{equation}\label{Eq:normf(x)}
 \partial f(x):=\tau_{\delta,\theta}(\partial x)\,,\quad\text{for each }x\in K.
 \end{equation}
Observe that the extension of~$\partial$ thus defined preserves all finite joins within each block~$K$ or~$f[K]$.

Pick $o\in K_\delta$ and fix a cofinal chain~$C$ in~$K$ with~$o$ as least element. We set
 \[
 F:=\setm{(c,u)\in C\times K}{u\utr c\text{ and }\delta\leq\partial u}
 \]
(cf. \eqref{Eq:defnutr}). Observe that the equality $c_{(\delta)}=u_{(\delta)}$ holds for each $(c,u)\in F$.

We refer the reader to Section~\ref{S:Basic} for the notation~$x^F$.

\setcounter{claim}{0}
\begin{claim}\label{Cl:FmsK2}
The subset $F$ is a cofinal \mss\ of~$K\times K$. Furthermore, the equality $(x,y)^F=(x^C\vee y^C,(x^C\vee y^C)_{(\delta\vee\partial y)})$ holds for each $(x,y)\in K\times K$.
\end{claim}

\begin{cproof}
As $(c,c)$ belongs to~$F$ for each $c\in C$, $F$ is cofinal in~$K\times K$. Let $(c,u),(d,v)\in F$, we must prove that $(c\wedge d,u\wedge v)\in F$. We may assume that $c\leq d$. It is straightforward to verify that $u\wedge v=c_{(\partial u)}\wedge d_{(\partial v)}=c_{(\partial u\wedge\partial v)}$, and so $(c\wedge d,u\wedge v)=(c,c_{(\partial u\wedge\partial v)})$ belongs to~$F$.

Now let $(x,y)\in K\times K$, put $\ol{c}:=x^C\vee y^C$ and $\eta:=\delta\vee\partial y$. As $o\vee y\leq\ol{c}$ and $\partial(o\vee y)=\eta$, it follows from Lemma~\ref{L:proj2meet}(i) that $\partial\ol{c}_{(\eta)}=\eta$, hence $(\ol{c},\ol{c}_{(\eta)})$ belongs to~$F$; it obviously lies above~$(x,y)$. For each $(c,u)\in F$ above~$(x,y)$, it follows from $x\leq c$, $y\leq u\leq c$, and $c\in C$ that $\ol{c}=x^C\vee y^C\leq c$. {}From $y\leq u$ it follows that $\partial y\leq\partial u$, thus, as $\delta\leq\partial u$, we get $\eta\leq\partial u$, and thus $\ol{c}_{(\eta)}\leq\ol{c}_{(\partial u)}\leq c_{(\partial u)}=u$ (the latter equality following from the relation $u\utr c$). Therefore, $(\ol{c},\ol{c}_{(\eta)})\leq(c,u)$.
\end{cproof}

As $F$ is a cofinal \mss\ of $K\times K$, it is a lower finite lattice; its least element is $(o,o)$.

\begin{claim}\label{Cl:F2ladd}
The lattice~$F$ is a $2$-ladder.
\end{claim}

\begin{cproof}
We prove that every element $(c,u)\in F$ has at most two lower covers in~$F$. We denote by~$c_*$ the unique lower cover of~$c$ in~$C$ if it exists (i.e., $c>o$), and by~$u^-$ the largest element of $\Proj(u)\setminus\set{u}$ such that $\delta\leq\partial u^-$ if it exists (i.e., $\partial u>\delta$). If~$c_*$ exists, then $(c_0,u_0):=(c_*,(c_*)_{(\partial u)})$ is an element of~$F$ smaller than~$(c,u)$. If~$u^-$ exists, then $(c_1,u_1):=(c,u^-)$ is an element of~$F$ smaller than~$(c,u)$.

Now let $(d,v)\in F$ such that $(d,v)<(c,u)$. If $d<c$, then $c_*$ exists and $d\leq c_*$, but $v\leq u$, thus
 \begin{align*}
 v&=d_{(\partial v)}&&(\text{because }v\in\Proj(d))\\
 &\leq(c_*)_{(\partial u)}&&
 (\text{because }d\leq c_*\text{ and }\partial v\leq\partial u)\,, 
 \end{align*}
and so $(d,v)\leq(c_0,u_0)$. Now suppose that $d=c$. Hence $v<u$, but $u,v\in\Proj(c)$, thus $v\in\Proj(u)$, so~$u^-$ exists and $v\leq u^-$, and so $(d,v)\leq(c_1,u_1)$. Therefore, each lower cover of~$(c,u)$ in~$F$ is equal to~$(c_i,u_i)$ for some~$i<2$.
\end{cproof}

Now we endow the disjoint union $\ol{K}:=(K\cup f[K])\sqcup F$ with the partial ordering obtained as the union of the respective partial orderings of the strong amalgam $K\cup f[K]$ and the $2$-ladder~$F$ together with the pairs
 \begin{align}
 x&<(c,u)\text{ just in case }x\leq c\,,\label{Eq:x<(cu)}\\
 f(x)&<(c,u)\text{ just in case }x\leq u\,,\label{Eq:f(x)<(cu)}
 \end{align}
for each $x\in K$ and each $(c,u)\in F$. It is important to observe that there is no ambiguity between~\eqref{Eq:x<(cu)} and~\eqref{Eq:f(x)<(cu)} in case $x\in I$ (i.e., $x=f(x)$), because of the equality $c_{(\delta)}=u_{(\delta)}$.
It is straightforward to verify that~$\ol{K}$ is a lower finite poset in which both~$K$ and~$f[K]$ are lower subsets. Furthermore, by using Claim~\ref{Cl:FmsK2}, we obtain that~$\ol{K}$ is a \jzs, in which the join operation on pairs of elements not in the same block~$K$ or~$f[K]$ is given by
 \begin{align}
 x_0\vee f(x_1)&=(x_0^C\vee x_1^C,(x_0^C\vee x_1^C)_{(\partial x_1)})\,,
 \label{Eq:x0veef(x1)}\\
 x\vee(c,u)&=(x^C\vee c,(x^C\vee c)_{(\partial u)})\,,\label{Eq:xveecy}\\
 f(x)\vee(c,u)&=(x^C\vee c,(x^C\vee c)_{(\partial x\vee\partial u)})\,,
 \label{Eq:fxveecy}\\
 (c_0,u_0)\vee(c_1,u_1)&
 =(c_0\vee c_1,(c_0\vee c_1)_{(\partial u_0\vee\partial u_1)})\,,
 \label{Eq:c01veey01}
 \end{align}
for all $x\in K$, all $x_0,x_1\in K\setminus I$, and all $(c,u),(c_0,u_0),(c_1,u_1)\in F$. For example, $x_0\vee f(x_1)$ is the least element of~$F$ above $(x_0,x_1)$, thus, by Claim~\ref{Cl:FmsK2} and as $\delta\leq\partial x_1$, it is equal to $(x_0^C\vee x_1^C,(x_0^C\vee x_1^C)_{(\partial x_1)})$. Observe that~\eqref{Eq:x0veef(x1)} implies immediately that every element of~$\ol{K}$ is the join of two elements of~$K\cup f[K]$---indeed, $(c,u)=c\vee f(u)$ for each $(c,u)\in F$.

We set
 \begin{equation}\label{Eq:normonF}
 \partial(c,u):=\tau_{\delta,\theta}(\partial u)\,,
 \quad\text{for each }(c,u)\in F\,.
 \end{equation}
In order to verify that the extension of~$\partial$ defined in~\eqref{Eq:normf(x)} and~\eqref{Eq:normonF} is a \jh, we consider the expressions~\eqref{Eq:x0veef(x1)}--\eqref{Eq:c01veey01} and we use repeatedly the fact that $\tau_{\delta,\theta}(\partial u)\geq\theta$ holds for each $(c,u)\in F$, hence~$\tau_{\delta,\theta}(\partial u)$ absorbs~$\partial x$ for each $x\in K$. The verifications of these facts are straightforward.

As $\partial\res_K$ has range~$\theta$, $\partial\res_{f[K]}$ has range  $\delta\cup(\theta'\setminus\theta)$, and~$\partial\res_F$ has range contained in (and in fact equal to)~$\theta'\setminus\theta$, we obtain that the range of~$\partial$ is~$\theta'$.

As both $K$ and~$f[K]$ are lower subsets of~$\ol{K}$ and $3$-ladders, in order to verify that~$\ol{K}$ is a $3$-ladder, it suffices to verify that each $(c,u)\in F$ has at most three lower covers in~$\ol{K}$. We use the notations~$c_*$ and~$u^-$ introduced in the proof of Claim~\ref{Cl:F2ladd}. By using Claim~\ref{Cl:F2ladd}, it is easy, although a bit tedious, to verify that every lower cover of~$(c,u)$ in~$\ol{K}$ has one of the following forms:
\begin{itemize}
\item $(c,u^-)$, $(c_*,(c_*)_{(\partial u)})$, and $f(u)$, if $c>o$ and $\partial u>\delta$;

\item $c$, $(c_*,(c_*)_{(\delta)})$, and $f(u)$, if $c>o$ and $\partial u=\delta$;

\item $o$ and $f(o)$, if $c=o$.
\end{itemize}
In any case, $(c,u)$ has at most three lower covers in~$\ol{K}$.
\end{proof}

\section{Getting a large $3$-ladder from a morass}\label{S:Mor2Ladd}

\begin{theorem}\label{T:Mor2Ladd}
Suppose that there exists an $(\omega_1,1)$-morass. Then there exists an atomistic $3$-ladder of cardinality~$\aleph_2$.
\end{theorem}

\begin{proof}
As discussed in Section~\ref{S:Morass}, there exists a simplified $(\omega_1,1)$-morass. We shall use the same notation as in Definition~\ref{D:Gap1Mor}. We shall construct inductively a system
 \[
 \cK=\biggl(\famm{K_\xi}{\xi\leq\omega_1},
 \famm{f^*}{f\in\bigcup\famm{\cF_{\alpha,\beta}}
 {\alpha<\beta\leq\omega_1}}\biggr)
 \]
satisfying the following conditions:
\begin{itemize}
\item[(K0)] Each $K_\xi$ is an atomistic normed $3$-ladder with range~$\theta_\xi$.

\emph{Accordingly, we shall denote by $\partial$, or $\partial_\xi$ in case~$\xi$ needs to be specified, the norm function on~$K_\xi$}.

\item[(K1)] $K_\xi$ is countable for each $\xi<\omega_1$.

\item[(K2)] If $\xi<\eta\leq\omega_1$ and $f\in\cF_{\xi,\eta}$, then $f^*\colon K_\xi\into K_\eta$ is a lower embedding.

\item[(K3)] If $\xi<\eta<\gamma\leq\omega_1$, then the equality $(f\circ g)^*=f^*\circ g^*$ holds for all $(f,g)\in\cF_{\eta,\gamma}\times\cF_{\xi,\eta}$.

\item[(K4)] If $\xi<\eta\leq\omega_1$, then the equality $\partial_\eta\circ f^*=f\circ\partial_\xi$ holds for each $f\in\cF_{\xi,\eta}$.
\end{itemize}

Provided we can carry out this construction we shall be done, since then $K_{\omega_1}$ is a normed $3$-ladder with range~$\omega_2$, thus of cardinality at least~$\aleph_2$. As every $3$-ladder has cardinality at most~$\aleph_2$ (cf. Proposition~\ref{P:UppBdBr}), $K_{\omega_1}$ has cardinality exactly~$\aleph_2$.

We start with $K_0:=\Pow(\set{0,1})$ (endowed with containment) with the norm operation defined by $\partial\es=0$, $\partial\set{0}=1$, and $\partial\set{1}=\partial\set{0,1}=2$.
Suppose the construction carried out to the level~$\alpha$ (that is, replace~$\leq\omega_1$ by~$\leq\alpha$ in (K0)--(K4) above), we shall show how to extend it to the level~$\alpha+1$. We put
 \begin{equation}\label{Eq:defIalpha}
 I_\alpha:=\setm{x\in K_\alpha}{\partial x<\delta_\alpha}\,.
 \end{equation}
By Lemma~\ref{L:amalg}, there exists a countable, transitive, normed $3$-ladder $K_{\alpha+1}$ with range $\theta_\alpha+(\theta_\alpha-\delta_\alpha)=\theta_{\alpha+1}$, containing~$K_\alpha$ as an ideal, together with a lower embedding $f_\alpha^*\colon K_\alpha\into K_{\alpha+1}$ such that the following conditions hold:
\begin{enumerate}
\item $K_\alpha\cup f_\alpha^*[K_\alpha]$ is the strong amalgam of~$K_\alpha$ and $f_\alpha^*[K_\alpha]$ over~$I_\alpha$;

\item every element of $K_{\alpha+1}$ is the join of two elements of $K_\alpha\cup f_\alpha^*[K_\alpha]$;

\item $f_\alpha^*\res_{I_\alpha}$ is the identity map on~$I_\alpha$;

\item The equality $\partial f_\alpha^*(x)=f_\alpha(\partial x)$ holds for each $x\in K_\alpha$ (\emph{recall that~$f_\alpha$ is the restriction of $\tau_{\delta_\alpha,\theta_\alpha}$ from~$\theta_\alpha$ into~$\theta_{\alpha+1}$}).
\end{enumerate}
In particular, it follows from the induction hypothesis and from (ii) that~$K_{\alpha+1}$ is also atomistic.

We put $\id_{\theta_\alpha}^*:=\id_{K_\alpha}$. Hence~$f^*$ is a lower embedding from~$K_\alpha$ into~$K_{\alpha+1}$, for each $f\in\cF_{\alpha,\alpha+1}$. It follows immediately from the definition that
 \begin{equation}\label{Eq:1step*eq}
 \partial_{\alpha+1}\circ f^*=f\circ\partial_\alpha\,,\quad
 \text{for each }f\in\cF_{\alpha,\alpha+1}\,.
 \end{equation}
We must define $f^*$, for each $\xi<\alpha$ and each $f\in\cF_{\xi,\alpha+1}$. It follows from~(P2) (cf. Definition~\ref{D:Gap1Mor}) that there are $g\in\cF_{\alpha,\alpha+1}$ and $h\in\cF_{\xi,\alpha}$ such that $f=g\circ h$. We have no other choice than to define $f^*:=g^*\circ h^*$.
It follows easily from Lemma~\ref{L:Mor2Tree} that the map~$h$ is uniquely determined, however the map~$g$ is not necessarily uniquely determined. This uniqueness problem arises if{f} $f=h=f_\alpha\circ h$, in which case the range of~$h$ is contained in the intersection of~$\theta_\alpha$ and the range of~$f_\alpha$, which is~$\delta_\alpha$. In this case, we need to prove that $f_\alpha^*\circ h^*=f_\alpha^*$, that is, $\rng h^*\subseteq I_\alpha$. This is obvious by~\eqref{Eq:defIalpha}, as $\rng h\subseteq\delta_\alpha$, so, by using the induction hypothesis, the relations
 \[
 \partial h^*(x)=h(\partial x)<\delta_\alpha
 \]
hold for each $x\in K_\xi$.

At this stage we are able to define $f^*=g^*\circ h^*$ as above. Once this is done, the verification of~(K3) is straightforward. This concludes the successor case.

Now let~$\lambda\leq\omega_1$ be a limit ordinal and suppose that the construction of~$\cK$ has been carried out to the level~$<\lambda$ (i.e., replace $\leq\omega_1$ by~$<\lambda$ in the formulations of (K0)--(K4)). The proof of this case runs as the proof of the limit case of the gap-$2$ transfer Theorem using a simplified morass presented in Devlin \cite[Section~VIII.4]{Devl84}, so we shall merely outline it. We put
 \[
 \cF_\lambda:=\bigcup\famm{\cF_{\alpha,\lambda}}{\alpha<\lambda}
 \qquad\text{(disjoint union)}\,,
 \]
and we define $d(f)$, for $f\in\cF_\lambda$, as the unique $\alpha<\lambda$ such that $f\in\cF_{\alpha,\lambda}$. For $f,f'\in\cF_\lambda$, let $f\tr f'$ hold, if $d(f)<d(f')$ and there exists a (necessarily unique) $g\in\cF_{d(f),d(f')}$ such that $f=f'\circ g$. It follows from the induction hypothesis that the map $\pi_{f,f'}=g^*\colon K_{d(f)}\into K_{d(f')}$ is a lower embedding and
 \begin{equation}\label{Eq:commpartial}
 \partial_{d(f')}\circ\pi_{f,f'}=g\circ\partial_{d(f)}\,.
 \end{equation}
It follows from~(P2) and~(P4) that~$\tr$ is an upward directed strict ordering on~$\cF_\lambda$, so we can form the direct limit
 \[
 \famm{K_\lambda,f^*}{f\in\cF_\lambda}=\varinjlim
 \famm{K_{d(f)},\pi_{f,f'}}{f\tr f'\text{ in }\cF_\lambda}
 \]
in the category of lattices and lattice homomorphisms. As all~$\pi_{f,f'}$s are lower embeddings, so are all $f^*$s. As~$K_\lambda$ is the directed union of all the ranges of the lower embeddings $f^*\colon K_{d(f)}\into K_\lambda$ and all the~$K_{d(f)}$s are atomistic $3$-ladders, $K_\lambda$ is an atomistic $3$-ladder. Furthermore, it follows from~\eqref{Eq:commpartial} that there exists a unique map~$\partial_\lambda\colon K_\lambda\to\theta_\lambda$ such that $\partial_\lambda\circ f^*=f\circ\partial_{d(f)}$ holds for each $f\in\cF_\lambda$, and it is straightforward to verify that~$\partial_\lambda$ is a \jh. As~$\partial_{d(f)}$ has range~$\theta_{d(f)}$ for each~$f\in\cF_\lambda$ and by~(P5), the map~$\partial_\lambda$ has range~$\theta_\lambda$. This completes the construction.
\end{proof}

As reminded in Section~\ref{S:Morass}, the nonexistence of a $(\omega_1,1)$-morass implies that~$\omega_2$ is inaccessible in the constructible universe~$\xL$. Hence,

\begin{corollary}\label{C:Mor2Ladd}
If there is no $3$-ladder of cardinality~$\aleph_2$, then~$\omega_2$ is inaccessible in~$\xL$.
\end{corollary}

We still do not know whether the nonexistence of a $3$-ladder of cardinality~$\aleph_2$ is consistent with~$\ZFC$. We thus reach the intriguing conclusion that the existence of a $3$-ladder of cardinality~$\aleph_2$ follows from either one of two axioms of set theory, namely $\MAP$ and the existence of a gap-$1$ morass, which do not imply each other. Indeed, as there are morasses in~$\xL$ and~$\xL$ satisfies the Generalized Continuum Hypothesis, while on the other hand $\MAP$ implies that $2^{\aleph_0}=2^{\aleph_1}$ (for the Martin-Solovay forcing, constructed from~$\aleph_1$ almost disjoint subsets of~$\omega$, used in the proof of Jech \cite[Theorem~16.20]{Jech03}, has precaliber~$\aleph_1$), the existence of a $(\omega_1,1)$-morass does not imply $\MAP$. A simpler argument, pointed by the referee, runs as follows. If $\mathsf{MA}(\aleph_1;\text{countable})$ (i.e., $\MA(\aleph_1)$ restricted to \emph{countable} notions of forcing) holds, then, following the argument of Jech \cite[Exercise~16.10]{Jech03} (using Cohen forcing), it is easy to prove that $2^{\aleph_0}>\aleph_1$. In particular, the existence of a gap-$1$ morass does not even imply $\mathsf{MA}(\aleph_1;\text{countable})$.

Conversely, Devlin~\cite{Devl78} finds a model of $\ZFC+\MA(\aleph_1)$ without any Kurepa tree, thus \emph{a fortiori} without any $(\omega_1,1)$-morass.

The question whether there exists a $k$-ladder of cardinality~$\aleph_{k-1}$, for $k\geq3$, is also raised in~\cite{Dito84}. This question seems to get harder as~$k$ grows larger (we have no proof of this, although the similar implication with lower finite lattices of breadth~$k$ is easy to establish). The proof of Theorem~\ref{T:Mor2Ladd} suggests (without proving it) that the existence of a gap-$n$ morass would entail the existence of a $(n+2)$-ladder of cardinality~$\aleph_{n+1}$, for every positive integer~$n$. However, such a program would run into a number of technical difficulties, starting with the lack of a formally established equivalence between gap-$n$ morasses and their simplified versions for $n\geq3$ (cf. Morgan~\cite{Morg98} for the case $n=2$), or even the lack of formally established proofs that the latter exist in~$\xL[A]$ for any $A\subseteq\omega_1$. Furthermore, our other result that $\MAP$ implies the existence of a $3$-ladder of cardinality~$\aleph_2$ (Theorem~\ref{T:MAP23ladd}) suggests that morasses may not even be the optimal tool required to solve that problem. Another objection is then that unlike morasses, higher and higher versions of Martin's Axiom would apparently not be able to construct ladders of cardinality beyond the continuum\dots\ Still, without the question whether $\ZFC$ implies the existence of a $3$-ladder of cardinality~$\aleph_2$ being settled, the morass track and the Martin's Axiom track seem to be all we have.

\section{Acknowledgment}
I am grateful to Rich Laver for inspiring early e-mail discussions about the consistency problem of $3$-ladders of cardinality~$\aleph_2$. I also thank Michael Pinsker for his many remarks and corrections, that lead in particular to substantial simplifications of several proofs.

\end{document}